\documentstyle[12pt]{article}
\title{On the concentration and the convergence rate with a moment condition in first passage percolation
} 
\author{Yu Zhang\footnote{Research supported by NSF grant DMS-4540247. }
\\
Department of Mathematics, University of Colorado}
\date{}
\begin{document}
\baselineskip .20in
\maketitle

\begin{abstract}
We consider the first passage percolation model on the ${\bf Z}^d$ lattice.
In this model, we assign  independently to each edge $e$ a non-negative passage time $t(e)$
with a common distribution $F$. Let $a_{0,n}$ be the passage time from the origin  to $(n,0,\cdots , 0)$.
Under the exponential tail assumption,
 Kesten (1993) and Talagrand (1995) investigated  the  concentration of $a_{0,n}$ from its mean using different methods.
With this concentration and the exponential tail assumption,
Alexander gave an estimate for  the convergence rate for ${\bf E} a_{0,n}$.
In this paper, focusing on  a moment condition, we reinvestigate the concentration  and  the convergence rate
for $a_{0,n}$ using a special  martingale structure.
\end{abstract}
\section {Introduction of the model and results.}
The first passage percolation model was introduced  in 1965 by Hammersley and  Welsh. 
In this model, we consider the ${\bf Z}^d$ lattice as a graph with edges connecting each {\em adjacent} or ${\bf Z}^d$-adjacent
pair of vertices
$u=(u_1, u_2,\cdots, u_d)$ and $v=(v_1, v_2,\cdots, v_d)$ 
with $d(u,v)=1$, where $d(u,v)$ is the Euclidean distance between $u$ and $v$.
For any two sets ${\bf A}$ and ${\bf B}$, we also define
$$d({\bf A}, {\bf B})= \min\{ d(u, v): u\in {\bf A}, v\in {\bf B}\}.$$
We assign independently to each edge a non-negative {\em passage time} $t(e)$ with a common distribution $F$.
More formally, we consider the following probability space. As the
sample space, we take $\Omega=\prod_{e\in {\bf Z}^d} [0,\infty),$ whose points 
are called {\em configurations}.
Let ${\bf P}=\prod_{e\in {\bf Z}^d} \mu_e$ be the corresponding product measure on $\Omega$, where
$\mu_e$ is the measure on $[0, \infty)$ such that
$$\mu_e(t(e)\leq x)= F(x). $$
The
expectation and variance with respect to ${\bf P}$ are denoted by ${\bf E}(\cdot)$
and $\sigma^2(\cdot)$, respectively.
For any two vertices $u$ and $v$, 
a path $\gamma$ from $u$ to $v$ is an alternating sequence 
$(v_0, e_1, v_1,...,v_i, e_{i+1}, v_{i+1},...,v_{n-1},e_n, v_n)$ of vertices $v_i$ and 
 edges $e_i$ between $v_i$ and $v_{i+1}$ in ${\bf Z}^d$ with $v_0=u$ and $ v_n=v$. 
Given such a path $\gamma$, we define its passage time  as 
$$T(\gamma )= \sum_{i=1}^{n} t(e_i).\eqno{(1.1)}$$
For any two sets ${\bf A}$ and ${\bf B}$, we define the passage time from ${\bf A}$ to ${\bf B}$
as
$$T({\bf A}, {\bf B})=\inf \{ T(\gamma)\},$$
where the infimum is over all possible finite paths from some vertex in ${\bf A}$ to some vertex in ${\bf B}$.
A path $\gamma$ from ${\bf A}$ to ${\bf B}$ with $T(\gamma)=T({\bf A}, {\bf B})$ is called the {\em optimal path} of 
$T({\bf A}, {\bf B})$. The existence of such an optimal path has been proven (see Kesten (1986) and Zhang (1995)) for $F(0)\neq p_c$, where $p_c$ is the critical probability for Bernoulli (bond) percolation on ${\bf Z}^d$.
We also want to point out that the optimal path may  not be  unique.
If we focus on a special configuration $\omega$, we may write $T({\bf A},{\bf B})(\omega)$
instead of $T({\bf A}, {\bf B})$. 
When ${\bf A}=\{u\}$ and ${\bf B}=\{v\}$ are single vertex sets, $T(u,v)$ is the passage time
from $u$ to $v$. We may extend the passage time over ${\bf R}^d$.
If $x$ and $y$ are in ${\bf R}^2$, we define $T(x, y)=T(x', y')$, where
$x'$ (resp., $y'$) is the nearest neighbor of $x$ (resp., $y$) in
${\bf Z}^2$. Possible indetermination can be eliminated by choosing an
order on the vertices of ${\bf Z}^2$ and taking the smallest nearest
neighbor for this order.

Given a vector $u\in {\bf R}^d$, by a subadditive argument, if ${\bf E}t(e) < \infty$, then
$$
\lim_{n\rightarrow \infty}{1\over n} T({\bf 0}, nu) = \inf_{n} {1\over n} {\bf E} T({\bf 0}, nu)=\lim_{n\rightarrow \infty}
{1\over n} {\bf E} T({\bf 0}, nu)=\mu_F(u) \mbox{ a.s. and in } L_1,\eqno{(1.2)}
$$
where the non-random constant $\mu_F(u)$ is  called the {\em time constant}. 
Later, Kesten showed (see Theorem 6.1 in Kesten (1986)) that
$$\mu_F(u)=0\mbox{ iff }F(0) \geq p_c.$$
In particular, Hammersley and  Welsh (1965) first studied the passage time $T({\bf 0}, nu)$ when $u=(1,0,\cdots, 0)$,
and  defined
$$a_{0,n}=T({\bf 0}, nu) \mbox{ and } \mu_F(u)=\mu_F.$$
When $F(0) < p_c$, the map $x \rightarrow \mu(x)$ induces  a norm on
${\bf R}^d$. The unit radius ball for this norm is denoted by ${\bf B}:={\bf B}(F)$
and is called the {\em asymptotic shape}. The boundary of ${\bf B}$ is
$$\partial {\bf B}:= \{ x \in {\bf R}^d: \mu(x)=1\}.$$
${\bf B}$ is a compact convex
deterministic set, and  $\partial {\bf B}$  is a continuous convex closed
curve (see Kesten (1986)).  Define for all $t> 0$,
$$B(t):= \{v\in {\bf R}^d,  T( {\bf 0}, v) \leq t\}.$$
We call $B(t)$ a  {\em random shape} on ${\bf Z}^d$.
The shape theorem (see Cox and Durrett (1981)) is the well-known result stating that for any
$\epsilon >0$,  if ${\bf E} t(e) < \infty$, then
$$t{\bf B}(1-\epsilon)  \subset {B(t) } \subset t{\bf B}(1+\epsilon ),
\mbox{ eventually w.p.1.}\eqno{(1.3)}$$

The natural and most challenging aspect in this field (see Kesten (1986) and  Smythe and Wierman (1978))
 is to question the ``speed" and  ``roughness"  of the  interface $B(t)$ 
from the deterministic boundaries $t{\bf B}$.  
Now we focus on the most interesting situation: when $F(0) < p_c$.
It is widely conjectured  (Kesten (1993)) that if
 $F(0) < p_c$, then
$$\sigma^2 (a_{0,n})\approx n^{2/3}.$$
The mathematical estimates for the upper bound of $\sigma^2(a_{0,n})$ are quite promising. Kesten (1993) 
showed that if $F(0)< p_c$ and ${\bf E}(t(e))^2 < \infty$, there is a constant $C_1$,
$$\sigma^2( a_{0,n})\leq C_1 n.\eqno{(1.4)}$$
In this paper,  $C$ and $C_i$ are always positive  constants that may depend on $F$, but not on $n$. 
Their values 
are not significant and  change from 
appearance to appearance.
Benjamini, Kalai, and  Schramm (2003) also showed that when $t(e)$ only takes two values $0<a< b$ with
a half probability for each one,
$$\sigma^2( a_{0,n})\leq C_1 n/\log n,\eqno{}$$
where $\log$ denotes the natural logarithm.
The probability estimate for the concentration has also been  investigated.
With an exponential  tail assumption:
$$\int^{\infty}_0 e^{\lambda x} dF(x) < \infty \mbox{ for some } \lambda >0,\eqno{(1.5)}$$
Kesten (1993) showed that
$${\bf P}\left [ \left |{a_{0,n} -{\bf E}a_{0,n}} \right|\geq x\sqrt{n}\right] \leq C_1 \exp(-C_2 x)
\mbox{ for } x\leq C_3 n.\eqno{(1.6)}$$
Later, under (1.5), Talagrand (1995) used his  isoperimetric inequality to show a better estimate as follows:
$${\bf P}\left [ \left |{a_{0,n} -{\bf E}a_{0,n} } \right|\geq x\sqrt{n}\right] \leq C_1 \exp(-C_2 x^2)
\mbox{ for } x\leq C_3 n.\eqno{(1.7)}$$

With (1.6) and the exponential  tail assumption in (1.5),
 Alexander (1993) investigates the convergence rate to show that 
$$n {\mu}_F \leq {\bf E}a_{0,n}\leq n {\mu}_F + C\sqrt{n} \log n.\eqno{(1.8)}$$

Clearly, if we reduce the tail assumption (1.5), for example,  with the $m$-th moment condition, we may not obtain 
the strong concentration result in (1.6). 
By a standard i.i.d. concentration estimate, we might believe that for any positive integer $m$,
$${\bf E}[ a_{0,n}-{\bf E}a_{0,n}]^{2m} \leq C n^m.\eqno{(1.9)}$$
In this paper, we show that (1.9) indeed holds with a logarithm correction as the following  theorem:\\

{\bf Theorem 1.} {\em If $F(0) < p_c$ and ${\bf E}(t(e))^{2m} < \infty$ for $m \geq 1$, then for any vector $u$
there exists 
a constant $C=C(d, u, F,m)$  such that}
$${\bf E}[T({\bf 0}, nu)-{\bf E}T({\bf 0}, nu)]^{2m} \leq C n^m \log^{10m} n.$$

In the proof of Theorem 1, the moment condition is only needed  for the edges around the origin and $u$. 
More precisely, for any small $\delta >0$, let
$$D_n(v)=v+[-3^dM\log^{1+\delta} n,  3^dM\log^{1+\delta} n],\eqno{(1.10)}$$
where $M$ is a constant to be given precisely after (2.5).
Now we consider the passage time  ${T}(D_n({\bf 0}), D_n(nu))$.
Note that
$${T}(D_n({\bf 0}), D_n(nu))\leq T({\bf 0}, nu) \leq \sum_{e\in D_n({\bf 0})\cup D_n(nu)} t(e)+ {T}(D_n({\bf 0}), D_n(nu)).\eqno{(1.11)}$$
By (1.11) and (1.2), if ${\bf E}t(e) < \infty$, then
$$
\lim_{n\rightarrow \infty}{1\over n} {T}(D_n({\bf 0}), D_n(nu))
=\mu_F(u) \mbox{ a.s. and in } L_1.\eqno{(1.12)}
$$
In this paper, we always assume that $d({\bf 0}, u)\leq 1$ without loss of generality when we work on ${T}(D_n({\bf 0}), D_n(nu))$.
With a weaker moment condition, we can show that
${T}(D_n({\bf 0}), D_n(nu))$ has a concentration property, as in the following theorem:\\

{\bf Theorem 2.}  {\em If $F(0) < p_c$ and ${\bf E}(t(e))^{1+\eta} < \infty$ for any $\eta >0$, then , for  $m\geq 1$,
there exists $C=C(d, \eta, M,m, F)$  such that 
$${\bf P}\left [ \left |{T}(D_n({\bf 0}), D_n(nu)) -
{\bf E}{T}(D_n({\bf 0}), D_n(nu)) 
\right|\geq \sqrt{n} \log^{7} n \right] \leq C_1 n^{-m}.\eqno{(1.13)}$$
Moreover, for all $m \geq 1$, there exists $C=C(F, d, \eta, M, m)$}
$$ {\bf E}\left( {T}(D_n({\bf 0}), D_n(nu)) -
{\bf E}{T}(D_n({\bf 0}), D_n(nu))\right)^{2m}\leq 
C n^m \log^{7m} n.$$

If we use Theorem 2 to reinvestigate the convergence rate for $a_{0,n}$, we only use  a moment condition to show (1.8),
 as in the following theorem:  \\

{\bf Theorem 3.} {\em If $F(0) < p_c$ and ${\bf E}(t(e))^{1+\eta} < \infty$ for any $\eta >0$, 
then there exists $C=C(F, \eta)$
such that}
$$ n\mu_F \leq {\bf E} a_{0,n}\leq n\mu_F+ C n^{1/2} \log^{4} n.$$

{\bf Remark 1.} It is possible to reduce the powers of the logarithm in Theorems 1-3.
However, the proofs are more complicated, so it may not be worth trying.\\

{\bf Remark 2.} Hammersley and Welsh (1965) introduced the following
point-face passage time:\\
$$ b_{0,n}=T({\bf 0}, {\bf L}_n) \mbox{ for } {\bf L}_n=\{(x_1,\cdots, x_d): x_1=n\}.$$
It was shown by them that
$$\lim_{n\rightarrow \infty} {b_{0,n}\over n}=\mu_F \mbox{ a.s. and in } L_1.$$
The same proofs in the following sections can be carried out to 
show the same results of Theorems 1-3 for $b_{0,n}$.
For the face-face passage time or the passage time of a minimal cutset
in a large box, we might need a different martingale structure.
Thus, we would like to explore these results in a different paper.

\section{ Preliminaries.} 

In  section 2, we would like to introduce a few basic first percolation  results for $F(0) < p_c$. 
Since $F(0) <p_c$, we may take $0<\epsilon <1$ small such that
$$F(\epsilon) < p_c.\eqno{(2.1)}$$
For each edge $e$, if $t(e) < \epsilon$, we say it is an $\epsilon^-$-edge; otherwise, it is an $\epsilon^+$-edge.
Let ${\bf C}_{\epsilon}({\bf 0})$ be an $\epsilon^-$-cluster containing the origin. By (2.1) and a standard estimate
(see Theorem 5.4 in Grimmett (1999)), there exist $C_1$ and $C_2$ such that
$${\bf P}[ |{\bf C}_{\epsilon}({\bf 0})| \geq m ]\leq C_1 \exp(-C_2 m).\eqno{(2.2)}$$

Now  we need to deal with the edges with a large passage time.
It is clear that there might be many  edges whose passage times are large. However, since $t(e)$ has a first moment, 
the cluster of these edges  cannot be very large. More precisely, since  ${\bf E} (t(e))^{1+\eta}< \infty$
by Markov's inequality, we select $M$ such that
$${\bf P}[t(e) \geq  M] \leq {{\bf E}t(e)\over M }\rightarrow 0 \mbox{ as } M\rightarrow \infty.\eqno{(2.3)}$$
For each edge $e$, if $t(e) \leq M$, we say it is an $M^-$-edge, otherwise, it is an $M^+$-edge.
If we take $M$ large, the probability (see Grimmett (1999)) that 
there exists a large $M^+$-cluster  is exponentially small.
However, this is not enough for our purpose. In fact, we need  any two points on the cluster boundary to be connected by a path with a passage time less than $M$. To work on this result, we need a few definitions.
For vertices $u$ and $v$, we say that $u$ and $v$ are ${\bf L}^d$-adjacent if 
$$d(u,v) < 2.$$
In other words, given a vertex $v$,
besides the vertices that are vertically and horizontally adjacent
to $v$ (${\bf Z}^d$-adjacent), we also consider the vertices that
are diagonally adjacent to $v$. Note that if $u$ and $v$
are ${\bf Z}^d$-adjacent, then they are ${\bf L}^d$-adjacent.
A path $(u_0,u_1,\cdots, u_n)$ is said to be an ${\bf L}^d$-path 
if $u_{i-1}$ and $u_i$ are ${\bf L}^d$-adjacent.
With this path, we can define an ${\bf L}^d$-cluster.

A vertex of ${\bf Z}^d$ is  {\em open}  if one of its $2d$ ${\bf Z}^d$-adjacent edges is  larger than $M$;
otherwise it is {\em closed}.  Note that if $v$ is closed, then 
the passage time of all of its ${\bf Z}^d$-adjacent edges must be equal  or
less than
$M$. Let $p(M)$ be the probability that a vertex  is open.
We denote by ${\bf C}_M({\bf 0})$  the open ${\bf L}^d$-cluster that contains the origin. Let
$$\theta(M)= {\bf P}[|{\bf C}_M({\bf 0})|=\infty] \mbox{ and } p_T(M)=\sup \{p(M): {\bf E} |{\bf C}_M({\bf 0})|< \infty\}.\eqno{(2.4)}$$
By a standard  computation (see chapter 2 in Grimmett (1999)), we can show that
$$0 < p_T(M) < 1.$$
In addition, if $p(M)< p_T(M)$,
we adopt from the same proof of Theorem 6.1 in Grimmett (1999) to show that
$${\bf P}[ |{\bf C}_{M}({\bf 0})| \geq m ]\leq C_1 \exp(-C_2 m).\eqno{(2.5)}$$
Note that 
$$p(M)= 1- ({\bf P}[ t(e) < M])^{2d}\rightarrow 0 \mbox{ as } M \rightarrow \infty,$$
so, from now on, we always take $M>1$ large such that $p(M)$ is  less than $p_T(M)$.

Now we introduce the boundary vertices  of a ${\bf Z}^d$- or  ${\bf L}^d$-cluster.
We consider  ${\bf A} $ as a finite ${\bf Z}^d$-cluster or a finite ${\bf L}^d$-cluster. If ${\bf A} $ is a ${\bf Z}^d$-cluster,
$v $ is said to be a boundary vertex of ${\bf A}$ when there exists $u \not\in {\bf A}$ and $v$
is ${\bf Z}^d$ adjacent to $u$. 
We denote by $\partial {\bf A}$  all boundary vertices of ${\bf A}$.
 We also let $\partial_o{\bf A} $, an outside boundary,  be  all the vertices  that are ${\bf Z}^d$-adjacent to ${\bf A}$
but not in ${\bf A}$.
Moreover, we denote by $\partial_o^e{\bf A}$ the edges between $\partial {\bf A}$ and $\partial_o{\bf A}$.
If ${\bf A} $ is an ${\bf L}^d$-cluster, we may also consider the boundary vertices by replacing ${\bf Z}^d$-connections  with ${\bf L}^d$-connections.
In addition, we denote by $\Delta {\bf A}$ the outside boundary vertices of $ {\bf A}$ such that there 
exist ${\bf Z}^d$-paths without using ${\bf A}$ from these vertices to $\infty$.
$\Delta {\bf A}$ is said to be the {\em exterior boundary} of ${\bf A}$.
By Lemma 1 in Kesten (1986), the exterior  boundary vertices
are ${\bf Z}^d$-adjacent. Note that all the  edges from
 a closed vertex have passage times less than $M$. 
On the other hand, since each vertex in the exterior boundary of
$ {\bf A}$ is ${\bf L}^d$-adjacent to ${\bf A}$, there are at most 
$3^d |{\bf A}|$ vertices in $\Delta {\bf A}$. 
We summarize these observations in the following lemma:\\

{\bf Lemma 1.} {\em Let ${\bf A}$ be a finite  open ${\bf L}^d$-cluster. For any two vertices $x$ and $y$ in its exterior boundary,
there exists a ${\bf Z}^d$-path $\gamma$  from $x$ and $y$ such that }
$$|\gamma| \leq 3^d|{\bf A}| \mbox{ and } t(e) \leq M \mbox { for each $e\in \gamma $}.$$

Now we define another  passage time $\{\tau(e)\}$ on ${\bf Z}^d$ as follows. We select $0< \delta < 1$ with 
the $\delta$ in (1.10). First, if
 $\epsilon \leq t(e) \leq  M$, then $\tau(e)=t(e)$.
Next,  if $t(e)< \epsilon$, but $e$ is not connected to
an $\epsilon^-$-cluster with more than $\log^{1+\delta} n$ vertices, 
then  we set $\tau(e)=t(e)$. But if $t(e)< \epsilon$, and $e$ is  
connected to an $\epsilon^-$-cluster with more than $\log^{1+\delta} n$ vertices, then  we set $\tau(e)=1$. 
Finally,
if $t(e) > M$, and the two vertices of $e$ are ${\bf L}^d$-adjacent  to an open ${\bf L}^d$-cluster with more than $\log^{1+\delta} n$ vertices, then 
$\tau(e)= 1$; otherwise, 
$\tau(e)= t(e).$ We want to comment that, unlike to
$\{t(e)\}$, $\{\tau(e)\}$ may not  be an i.i.d. sequence.
To see this, suppose that the two vertices of $e$ are connected to two separated open clusters with a size of
$2\log^{1+\delta} (n)/3$. If $t(e)>M$, then the value of these edges $\{e'\}$ in these clusters will be $\tau(e')=1$.
But if $t(e) < M$, then these values will stay the same as $t(e')$.

For a ${\bf Z}^d$-path $\gamma$, we define 
$$T_\tau (\gamma)=\sum_{e\in \gamma} \tau(e).$$
For any two vertices sets  ${\bf A}$ and ${\bf B}$, we denote by
$$T_\tau({\bf A}, {\bf B}) =  \inf \{ T_\tau(\gamma): \gamma  \mbox{ a path from ${\bf A} $ to ${\bf B}$}\},\eqno{(2.6)}$$
where the path is a ${\bf Z}$-path. In fact, after Lemma 1,
we barely use the ${\bf L}^d$-connection. 
With this passage time, let
$$B_\tau(t)= \{v\in {\bf Z}^d: T_\tau(D_n({\bf 0}), v) \leq t\}.$$
Thus, for each $k\geq 0$,
$$D_n({\bf 0})\subseteq B_\tau(k) \subseteq B_\tau(k+1).\eqno{(2.7)}$$

We would like to introduce a few geometric properties for $B_\tau(k)$.
Note that the largest $\epsilon^-$-cluster is at most $\log^{1+\delta} n$,  so $B_\tau(k)$ is finite for each $k$.  
With the definition of $B_\tau(k)$, it is easy to verify    the following arguments in (2.8)--(2.9): 
$$B_\tau(k)\mbox{ is ${\bf Z}^d$-connected}.\eqno{(2.8)}$$
$$\mbox{ $T_\tau (D_n({\bf 0}), v)\leq k$ for $v\in B_\tau(k)$, and  
$T_\tau(D_n({\bf 0}), u)> k$ for 
$u\not\in B_\tau(k)$.}\eqno{(2.9)}$$
In particular, if $v\in \partial_o B_\tau(k)$, then 
$$T_\tau(D_n({\bf 0}), v)> k.\eqno{(2.10)}$$
Let $\gamma$ be an optimal path for $T_\tau(D_n({\bf 0}), D_n (nu))$.
Here we assume that $n$ is large enough such that
$$D_n({\bf 0})\cap D_n(nu)=\emptyset; $$
otherwise,
$$T_\tau(D_n({\bf 0}), D_n (nu))=0.$$
Now we assume that (see Fig. 1) 
$$D_n(nu)\cap  B_\tau(k)=\emptyset \mbox{ for some $k$}. $$
We go along an optimal path $\gamma$ from $D_n({\bf 0})$ to meet $\partial B_\tau(k)$ at $w$, and then from $w$ to $D_n(nu)$. 
Thus,
$$k+T_\tau(\partial B_\tau(k), D_n(nu)) \leq T_\tau(\gamma)=T_\tau(D_n({\bf 0}), D_n (nu)).\eqno{(2.12)}$$
On the other hand, we select an optimal path $\gamma'$
 from $\partial B_\tau(k)$ to $D_n(nu)$. 
Let $v$ be the intersection of the path and 
$\partial B_\tau(k)$. 
By (2.9),  there exists a path $\gamma''$ from $v$ to $D_n({\bf 0})$ (see Fig. 1)
with 
$$T(\gamma'') \leq k.$$
Thus, 
$$ T_\tau(D_n({\bf 0}), D_n (nu))\leq T_\tau(\gamma')+ T_\tau(\gamma'')\leq  k+T_\tau ( \partial B_\tau(k), D_n(nu)).\eqno{(2.13)}$$
Note that
$$T_\tau ( \partial B_\tau(k),D_n(nu))=T_\tau (B_\tau(k),D_n(nu)),\eqno{(2.14)}$$
so together with (2.12) and (2.13), we have the following lemma.\\

{\bf Lemma 2.} {\em If $D_n(nu)\cap  B_\tau(k)=\emptyset $,  then}
$$ k+T_\tau ( B_\tau(k), D_n(un))=  T_\tau(D_n({\bf 0}), D_n(nu)).$$

Now we will give an upper bound for the passage time
of each edge in an optimal path. Let $\gamma$ be an optimal
path from $D_n({\bf 0})$ to $ D_n(nu)$.
We select an edge $e$ from $\gamma$ with the two vertices
$v$ and $v'$ and let
$$\gamma=\gamma'' \cup \{e\} \cup \gamma' \mbox{ (see Fig. 1).}$$
If $\tau(e ) > 3^dM \log^{1+\delta} n$,  
then the open ${\bf L}^d$-cluster including  $v$ and $v'$ is
 smaller than $\log^{1+\delta} n$, so by
Lemma 1, we use pieces of $\gamma'$ and $\gamma''$ and the edges in the exterior boundary of $\Delta {\bf C}_M(v)$ (see Fig. 1) to 
construct another path from $D_n({\bf 0})$ to 
$D_n(nu)$.  Note that if $v$ is near $D_n({\bf 0})$
or $D_n(nu)$, we may not need $\gamma'$ or $\gamma''$.
By Lemma 1, the passage time of each edge in the path of
 the exterior boundary of $\Delta {\bf C}_M(v)$
is less than
$M$. By Lemma 1 again,
$$|\Delta {\bf C}_M(v)|\leq 3^d \log^{1+\delta} (n).$$
With these observations, 
the passage time in this newly constructed path is strictly less than the passage time of the original optimal path (see Fig. 1).
The contradiction tells us that 
$ \tau(e)\leq 3^d M \log^{1+\delta} n$.
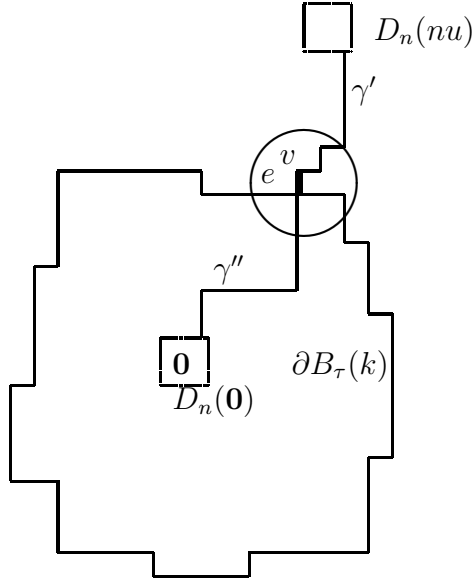
\begin{figure}\label{F:graphG}
\begin{center}
\setlength{\unitlength}{0.0125in}%
\begin{picture}(200,180)(67,910)
\thicklines
\put(130,995){\dashbox{.5}(20,20)[t]}
\put(190,1135){\dashbox{.5}(20,20)[t]}

\put(190,1080){\circle{1000}{$$}}
\put(187,1075){\line(0,1){10}}
\put(187.5,1075){\line(0,1){10}}
\put(188,1075){\line(0,1){10}}
\put(188.5,1075){\line(0,1){10}}
\put(189,1075){\line(0,1){10}}
\put(187,1085){\line(1,0){10}}
\put(197,1085){\line(0,1){10}}
\put(197,1095){\line(1,0){10}}
\put(207,1095){\line(0,1){40}}
\put(180,1088){$\small{v}$}
\put(172,1080){$\small{e_{}}$}
\put(210,1115){$\gamma'$}
\put(215,1140){$\mbox{ $D_n(nu)$}$}
\put(135,985){${{D_n({\bf 0})}}$}

\put(135,1000){${{\bf 0}}$}
\put(152,1040){{$\gamma''$}}
\put(187,1075){\line(0,-1){40}}
\put(187,1035){\line(-1,0){40}}
\put(147,1035){\line(0,-1){20}}

\put(187,1075){\line(1,0){20}}
\put(207,1075){\line(0,-1){20}}
\put(207,1055){\line(1,0){10}}
\put(217,1055){\line(0,-1){30}}
\put(217,1025){\line(1,0){10}}
\put(227,1025){\line(0,-1){60}}
\put(227,965){\line(-1,0){10}}
\put(217,965){\line(0,-1){40}}
\put(217,925){\line(-1,0){50}}
\put(167,925){\line(0,-1){10}}
\put(167,915){\line(-1,0){40}}
\put(127,915){\line(0,1){10}}
\put(127,925){\line(-1,0){40}}
\put(87,925){\line(0,1){30}}
\put(87,955){\line(-1,0){20}}
\put(67,955){\line(0,1){40}}
\put(67,995){\line(1,0){10}}
\put(77,995){\line(0,1){50}}
\put(77,1045){\line(1,0){10}}
\put(87,1045){\line(0,1){40}}
\put(87,1085){\line(1,0){60}}
\put(147,1085){\line(0,-1){10}}
\put(147,1075){\line(1,0){40}}
\put(185,1000){${\partial B_\tau(k)}$}

\end{picture}
\end{center}
\caption{\em The graph shows how to use $\gamma'$ and $\gamma''$  to construct
a path from $D_n({\bf 0})$ to $D_n(nu)$. This shows that
$ T_\tau(D_n({\bf 0}), D_n (nu))\leq T_\tau(\gamma')+ T_\tau(\gamma'')$.
The circuit is the boundary of the $M^+$-cluster of $v$. The number of edges in the circuit
is less than $3^d\log^{1+\delta} n$, and the passage time of each edge is less than $M$. We go along $\gamma''$ from $D_n({\bf 0})$, a dashed box, to reach the circuit, then along the circuit 
to reach $\gamma'$, and finally along $\gamma'$ to reach $D_n(nu)$, another dashed box. The passage time of this new constructed 
path is strictly less than the passage time of the original optimal path if $\tau(e) > M3^d \log^{1+\delta} n.$ }
\end{figure}

We summarize these observations  by the following lemma:\\

{\bf Lemma 3.} {\em 
For each edge $e$ in an optimal path of $T_\tau(D_n({\bf 0}), D_n(nu))$, }
$$\tau(e) \leq C\log^{1+\delta} n.\eqno{(2.15)}$$

Given a fixed connected set $\Gamma$ containing the origin, we define the 
 event as
$$\{B_\tau(k)=\Gamma\}=\{\omega: B_\tau(k)(\omega)=\Gamma\}.\eqno{(2.16)}$$
The following lemma is directly from the definition:\\

{\bf Lemma 4.} {\em If $\Gamma_1$ and $\Gamma_2$ are two different vertex sets, then}
$$\{B_\tau (k)=\Gamma_1\}\cap \{B_\tau (k)=\Gamma_2\}=\emptyset.\eqno{}$$

Since the $\epsilon^-$-cluster cannot be larger than $\log^{1+\delta} (n)$, 
any path from $D_n({\bf 0})$  to the boundary of 
$[- 2\log^{1+\delta} n,  2\log^{1+\delta} n]^d$ costs  at least passage time $1$.
Furthermore, any path from the origin to the boundary of 
$[-3 \log^{1+\delta} n,  3\log^{1+\delta} n]^d$ costs at least passage time $2$.
By a simple induction,  any path from the origin to the boundary of $[-k \log^{1+\delta} n,  k\log^{1+\delta} n]^d$
costs at least passage time $k-1$. With this observation, we have the following lemma:\\

{\bf Lemma 5.} {\em For each $k$,
$$B_\tau(k) \subset [- (k+1)\log^{1+\delta} n, (k+1)\log^{1+\delta} n]^d\mbox{ and }
d(B_\tau(k), B_\tau(k+1))\leq C\log^{1+\delta} n.\eqno{}$$}

In Lemma 5, we showed that $B_\tau(k)$ in a cube. The following 
lemma will show another direction.\\

{\bf Lemma 6.} {\em For each vector $u$, there exists $C=C(F, d, u)$
such that}
$$ B_\tau(C n\log^{1+\delta}(n))\cap D_n(nu)\neq \emptyset.$$

{\bf Proof.} There exists a path $\gamma$ (not necessarily optimal) from
$D_n({\bf 0})$ to $D_n(nu)$ with 
$$|\gamma|\leq Cn.$$
For each edge $e$ in $\gamma$, if $\tau(e) > M$, 
then we use the trick in (2.15) to
avoid using $e$ by
using at most $3^d \log^{1+\delta}(n)$ other edges.
Note that the passage time of these edges
 is less than or equal to $M$. The newly constructed path is denoted
 by $\gamma'$. If there exists an edge $e'$ in $\gamma'$ with 
$\tau(e') > M$, we use
the same argument to work on  $e'$. Thus, by a simple induction, we can find a path 
$\bar{\gamma}$ from $D_n({\bf 0})$ to $D_n(nu)$ such that
$$|\bar{\gamma}| \leq C n 3^d \log^{1+\delta}(n).$$
On the other hand, the passage time of each edge in $\bar{\gamma}$
is less than or equal to $M$.
Therefore, Lemma 6 follows. $\Box$\\

The following lemma shows  that optimal paths have to stay inside $B_\tau(k)$ if $D_n(nu)\cap  B_\tau(k)\neq \emptyset$:\\

{\bf Lemma 7.} {\em If  $D_n(nu)\,\,\cap  \,\,B_\tau(k)\neq \emptyset$ for some $k$, then all the optimal paths for $T_\tau(D_n({\bf 0}), D_n(nu))$
 stay inside $B_\tau(k)$.}\\

{\bf Proof.} Suppose that Lemma 7 does not hold. 
There exists
a path $\gamma$ for $T_\tau(D_n({\bf 0}), D_n(nu))$
such that some of its vertices in  $\partial_o B_\tau(k)$ and
$$T_\tau(\gamma) =T_\tau(D_n({\bf 0}), D_n(nu)).\eqno{(2.17)}$$
Note that $D_n(nu)\cap  B_\tau(k)\neq \emptyset$, so 
$$T_\tau(D_n({\bf 0}), D_n(nu))\leq k.\eqno{(2.18)}$$
However, by (2.9), since some of $\gamma$'s vertices are  in $\partial_o B_\tau(k)$,
$$T_\tau(\gamma) > k,$$
a contradiction. $\Box$\\

We will show the following lemma 
to demonstrate the relationship between $T_\tau(D_n({\bf 0}), D_n(nu))$ and ${T}(D_n({\bf 0}), D_n(nu))$:\\

{\bf Lemma 8.} {\em If $F(0) < p_c$ and $p(M) < p_T(M)$, then for all $u$ with $d({\bf 0}, u)\leq 1$, there exist
$C_i=C_i(F, M, \epsilon)$ for $i=1,2$ such that}
$${\bf P} [|{T}(D_n({\bf 0}), D_n(nu))-T_\tau(D_n({\bf 0}), D_n(nu)) |>0]\leq C_1 \exp(-C_2\log^{1+\delta}n).\eqno{(2.19)}$$

{\bf Proof.} By Proposition 5.8 in Kesten (1986), with a probability larger than $1-C_1\exp(-C_2 n)$,
there exists an optimal path $\gamma$ for ${T}(D_n({\bf 0}), D_n(nu))$ with $\gamma\leq L n$. 
On the existence of $\gamma$,
by the definition of $\tau$-edges,  if
$$|{T}(D_n({\bf 0}), D_n(nu))- T_\tau({\bf 0}), D_n(nu))|>0,$$ 
then there exists $v\in \gamma$ such that
$$|{\bf C}_{\epsilon^-}(v)|\geq \log^{1+\delta} n \mbox{ or } |{\bf C}_{M^+}(v)|\geq \log^{1+\delta} n.\eqno{(2.20)}$$
Thus, note that there are at most $(Ln)^{2d}$ choices for $v$,
so Lemma 8 follows from (2.3) and (2.5).  $\Box$\\

As we mentioned,  $\{\tau(e)\}$ is not an i.i.d.
sequence. 
Fortunately, when two edges $e_1$ and $e_2$ are separated far away, $\tau(e_1)$ and $\tau(e_2)$ are independent.\\

{\bf Lemma 9.}
{\em Let $v(e_1)$ and $v(e_2)$ be two vertices of $e_1$ and $e_2$ with
$$d(v(e_1), v(e_2))> 2\log^{1+\delta}(n). \eqno{(2.22)}$$
Then
$\tau(e_1)\mbox{ and } \tau(e_2) \mbox{ are independent}. $}\\

{\bf Proof.} If $\epsilon \leq t(e_1)\leq M$, then $\tau(e_1)$ keeps the same value regardless of $t(e')$ for
$e\neq e'$.
If $e_1$ is not connected to an $\epsilon^-$-cluster or an
$M^+$-cluster with a size larger than $\log^{1+\delta}(n)$, then  $\tau(e_1)=t(e_1)$ will also keep the same value
regardless of  $t(e')$ for
$$d(v(e_1), v(e'))> \log^{1+\delta}(n).$$ 
If $e_1$ is  connected to an $\epsilon^-$-cluster or
an $M^+$-cluster with a size larger than $\log^{1+\delta}(n)$, then $\tau(e_1)=1$.  Note that
the size of the cluster is always  larger than $\log^{1+\delta}(n)$ without edge $e'$ when
$$d(v(e_1), v(e'))> \log^{1+\delta}(n).$$ 
Therefore, $\tau(e_1)$ is always 1 regardless of  $t(e')$.
In summary, the values of $\tau(e_1)$ only depend
on the values of   $\{t(e')\}$ if 
$$d(v(e_1), v(e'))\leq \log^{1+\delta}(n).$$ 
The same argument also works for $\tau(e_2)$.
Since $\{t(e)\}$ is i.i.d.,  $\tau(e_1)$ and  $\tau(e_2)$ are independent if (2.22) holds. Lemma 9 follows. $\Box$\\

Given a vertex set $\Gamma$, we define its {\em shell}  as follows. Let (see Fig. 2)
$$S(\Gamma)=\{u\in {\bf Z}^d:  2\log^{1+\delta} (n) < d(u, \Gamma)< 1+2\log^{1+\delta} (n)\}.$$
We also define the sets outside the shell as follows (see Fig. 2):
$${S}^+(\Gamma)=\{u\in {\bf Z}^d:  2\log^{1+\delta} (n) < d(u, \Gamma)\}.$$
By using Lemma 9, we have the following lemma.\\
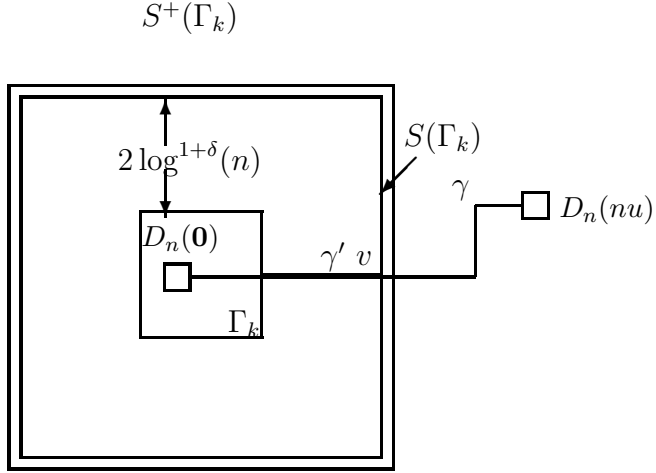
\begin{figure}\label{F:graphG}
\begin{center}
\setlength{\unitlength}{0.0125in}%
\begin{picture}(200,180)(67,910)
\thicklines
\put(150, 1000){\framebox(50,50)[br]{\mbox{$\Gamma_k$}}}
\put(100, 950){\framebox(150,150)[br]{\mbox{$$}}}
\put(95, 945){\framebox(160,160)[br]{\mbox{$$}}}
\put(140,1070){$2\log^{1+\delta}(n) $}
\put(160,1080){\vector(0,1){20}}
\put(160,1070){\vector(0,-1){20}}
\put(260,1080){$S(\Gamma_k)$}
\put(265,1075){\vector(-1,-1){15}}
\put(150,1130){$S^+(\Gamma_k)$}
\put(170, 1025){\line(1,0){120}}
\put(200, 1026){\line(1,0){50}}

\put(240, 1030){$v$}
\put(225, 1030){$\gamma'$}

\put(290, 1025){\line(0,1){30}}
\put(290, 1055){\line(1,0){20}}
\put(310, 1050){\framebox(10,10)[br]{\mbox{$$}}}
\put(325, 1050){\small{$D_n({nu})$}}
\put(280, 1060){$\gamma$}

\put(160, 1020){\framebox(10,10)[br]{\mbox{$$}}}
\put(150, 1038){\small{$D_n({\bf 0})$}}

\end{picture}
\end{center}
\caption{\em The graph shows the shell of $\Gamma_k$ and the vertex set outside of the shell. There is an optimal path 
$\gamma$ from  $D_n(nu)$ to the shell $S(\Gamma_k)$ at $v$. The path $\gamma'$ is from $v$ to $\Gamma_k$ with a length
less than $2\log^{1+\delta}(n)+1$. }
\end{figure} 

{\bf Lemma 10.}  {\em Let $\Gamma_k$ be vertex set containing the origin. If $\kappa\subset  {S}^+(\Gamma_k)$, then
$ T_\tau(S(\Gamma_k), \kappa )$ and $I( B_\tau(k)=\Gamma_k)$ are independent, where $I({\cal A})$ is the indicator of ${\cal A}$.}\\

{\bf Proof.} 
It follows from the same proof of Zhang's Proposition  3 (2005) that the event of
$\{B_\tau(k)=\Gamma_k\}$ only depends on the configurations of $\tau(e)$ for $e\in \Gamma_k \cup \partial_o \Gamma_k$. On the other hand, $T_\gamma(S(\Gamma_k), \kappa)$ only depends on the configurations of 
$\tau(e)$ for $e\in S(\Gamma_k)$. By Lemma 9,
$\tau(e)$ for $e\in \Gamma_k\cup \partial_o\Gamma_k$ and $\tau(e')$ for $e\in S^+(\Gamma_k)$ are 
independent. Therefore, Lemma 10 follows. $\Box$

\section{ Constructing a martingale sequence.}
Kesten (1993) constructed  a martingale to investigate
the concentration for passage time. Later,
 Kesten and Zhang (1997) and  Higuchi and Zhang (2000) used
 different martingale structures to investigate
the concentration for passage time. In this paper,
we use the sets $\{B_\tau(k)\}$ to construct a martingale
sequence.

Let 
$${\cal F}_{-1}=\{\emptyset, \Omega\} \mbox{ be the trivial $\sigma$-field}.$$
Let 
$\bar{B}_\tau(k)$ be all the edges with vertices in
 ${B}_\tau(k)\cap \partial_o B_\tau(k)$. 
With this edge set, let 
$${\cal F}_k\mbox{  be the $\sigma$-field generated by }\{\tau(e) : e \in \bar{B}_\tau(k)\}. \eqno{(3.1)}$$
More precisely, ${\cal F}_k$ is the smallest $\sigma$-field that contains all the sets  with the form  
$$\{\omega\in \Omega: \bar{B}_\tau(k)(\omega)=\Gamma_k, (\tau(e_1)(\omega),\cdots, \tau(e_i)(\omega))\in {\bf B}\}, $$
where  $\Gamma_k=(e_1,\cdots, e_i)$ for edges $\{e_j\}$ $j=1,2\cdots, i$, and ${\bf B}$ is an $i$-dimensional Borel set. Note that 
$B_\tau(k) \subseteq  B_\tau(k+1)$, so 
$${\cal F}_k \subseteq {\cal F}_{k+1}.$$

With this filtration,
$${\cal F}_{-1}\subset {\cal F}_0\subset {\cal F}_1\subset \cdots \subset{\cal F}_k\subset \cdots,$$
we set the following  martingale sequence:
$$\{M_k={\bf E}(T_\tau(D_n({\bf 0}), D_n(nu))\,\, |\,\,{\cal F}_k)- {\bf E}T_\tau(D_n({\bf 0}), D_n(nu))\}\mbox{ for } k=-1,0,1,2,\cdots .$$
The martingale differences of $M_k$ are denoted by
$$\Delta_{k,n}=M_{k+1}-M_k= {\bf E}(T_\tau(D_n({\bf 0}), D_n(nu))\:\:|\:\:{\cal F}_{k+1})-{\bf E}(T_\tau(D_n({\bf 0}), D_n(nu))\:\:|\:\:{\cal F}_k).$$
The following  work in  section 3 estimates the martingale difference.
Since $M_{-1}=0$, we estimate the martingale differences when 
$k \geq 0$.

Now we divide the shape of $B_\tau(k)$ into a few situations to estimate $|\Delta_{k,n}|$. 
First, we assume that 
$$D_n(nu)\cap {B}_\tau(k)\neq \emptyset .$$
 Let $I_1(k)$ be the indicator of this event. Note that
 $$\{ D_n(nu)\cap {B}_\tau(k)\neq \emptyset \} =\{ \exists
 \,\,\, \gamma\subset B_\tau(k): T_\tau(\gamma)\leq k\}.$$
 In other words, the event only depends on the configurations
 of edges inside $B_\tau(k).$ Thus, $I_1(k)$ is ${\cal F}_k$-measurable.
By Lemma 7, all optimal paths will stay inside $B_\tau(k)$. Thus, 
$T_\tau(D_n({\bf 0}), D_n(nu))I_1(k)$ is ${\cal F}_l$-measurable for all $l\geq k$.
Therefore,  
$$\Delta_{k,n}= {\bf E}(T_\tau(D_n({\bf 0}), D_n(nu))I_1(k)\:\:|\:\:{\cal F}_{k+1})-{\bf E}(T_\tau(D_n({\bf 0}), D_n(nu))I_1(k)\:\:|\:\:{\cal F}_k)=0.\eqno{(3.2)} $$
By  Lemma 6,  there exists a constant $C$ such that
$$D_n(nu)\cap {B}_\tau(C n \log^{1+\delta}(n))\neq \emptyset .$$ Thus, by Lemma 7, 
all the  optimal paths of $T_\tau(D_n({\bf 0}), D_n(nu))$
will stay inside $B_\tau(Cn\log^{1+\delta}(n))$. 
For  all $u$ with $d({\bf 0}, u) \leq 1$ and $n$, if $l \geq Cn\log^{1+\delta}(n)$, the same argument of (3.2) implies that
$$\Delta_{k,n}= [{\bf E}(T_\tau(D_n({\bf 0}), D_n(nu))\:\:|\:\:{\cal F}_{l+1})-{\bf E}(T_\tau(D_n({\bf 0}), D_n(nu))\:\:|\:\:{\cal F}_l)]=0.\eqno{(3.3)} $$

Now we deal with the most difficult situation, 
$$\{D_n(nu)\cap {B}_\tau(k)= \emptyset\}.\eqno{(3.4)}$$ 
Let $I_2(k)$ be the indicator of this event. 
Note that $I_1(k)$ is ${\cal F}_k$-measurable and 
$$I_2(k)=1-I_1(k),\eqno{(3.5)}$$
so
$$I_2(k) \mbox{ is ${\cal F}_k$-measurable, and, thus, is ${\cal F}_{k+1}$-measurable}.\eqno{(3.6)}$$
Note that
$$T_\tau({B}_{\tau}(k), D_n(nu )) (1-I_2(k))=0,$$
so by Lemma 2 and (3.6)
\begin{eqnarray*}
&&kI_2(k)+{\bf E} (T_\tau({B}_{\tau}(k), D_n(nu ))\,\,| \,\,{\cal F}_{k+1})\\
&=&kI_2(k)+{\bf E} (T_\tau({B}_{\tau}(k), D_n(nu ))I_2(k)\,\,| \,\,{\cal F}_{k+1})\\
&= &{\bf E}[T_\tau(D_n({\bf 0}), D_n(nu))I_2(k) \,\,| \,\,{\cal F}_{k+1}].\hskip 8.6cm (3.7)
\end{eqnarray*}
Similarly, 
\begin{eqnarray*}
&&kI_2(k)+{\bf E} (T_\tau({B}_{\tau}, D_n(nu ))\,\,| \,\,{\cal F}_{k})\\
&=&kI_2(k)+{\bf E} (T_\tau({B}_{\tau}(k), D_n(nu ))I_2(k)\,\,| \,\,{\cal F}_{k})\\
&=&{\bf E}[T_\tau(D_n({\bf 0}), D_n(nu))I_2(k) \,\,| \,\,{\cal F}_{k}].\hskip 8.7cm (3.8)
\end{eqnarray*}

By (3.7) and (3.8), 
\begin{eqnarray*}
&&|{\bf E}[ T_\tau(D_n({\bf 0}), D_n(nu))I_2(k)\,\,|\,\, {\cal F}_{k+1}]
 -{\bf E}[T_\tau(D_n({\bf 0}), D_n(nu))I_2(k)\,\,|\,\,{\cal F}_k]|\\
&= & |{\bf E} [T_\tau(B_\tau(k), D_n(nu) )\,\,| \,\,{\cal F}_{k+1}]-
{\bf E} [T_\tau(B_\tau(k), D_n(nu))\,\,| \,\,{\cal F}_k]|.
\hskip 4cm (3.9)
\end{eqnarray*}

Now we estimate 
${\bf E} [T_\tau(B_\tau(k), D_n(nu))\,\,| \,\,{\cal F}_k]$. 
For fixed vertex sets 
${\bf 0} \in \Gamma_k$,
 let 
$I(\Gamma_k )$ be the indicator of the event that  $\{B_\tau(k)= \Gamma_k\}.$
Thus, by Lemma 4,
$$
{\bf E} [T_\tau(B_\tau(l), D_n(nu) )\,\,| \,\,{\cal F}_k]=
\sum_{\Gamma_l}{\bf E} [T_\tau( \Gamma_l, D_n(nu) )I(\Gamma_k )\,\,| \,\,{\cal F}_k)],\eqno{(3.10)}
$$
where the sum takes all possible vertex sets containing the origin.
Using  the definition of the shell in Lemma 8, we have,
$$ T_\tau(S(\Gamma_{k}),  D_n(nu))\leq T_\tau(\Gamma_{k},  D_n(nu)).\eqno{(3.11)}$$
If $S(\Gamma_{k})\cap D_n(nu)=\emptyset$,
let $\gamma$ be an optimal path starting at $v$ from $S(\Gamma_{k})$ to $D_n(nu)$ (see Fig. 2). 
On the other hand, note that 
$$d(\Gamma_{k}, S(\Gamma_{k})) \leq C\log^{1+\delta} (n),$$
so there exists a path $\gamma'$ from $v$ to $\Gamma_{k}$ with less than $C\log^{1+\delta} (n)$ edges (see Fig. 2).
For each edge $e\in \gamma'$, if  $\tau(e)$ is larger than $M$, we may use the rerouting argument in Lemma 3
with extra $C_1\log^{1+\delta}(n)$ passage 
time to avoid using $e$. With this observation, we can show that
$$  T_\tau( \Gamma_{k},  D_n(nu))\leq T_\tau(S(\Gamma_{k}),  D_n(nu))+C \log^{2+2\delta}(n).\eqno{(3.12)}$$
If $S(\Gamma_{k})\cap D_n(nu)\neq \emptyset$, then $T_\tau(S(\Gamma_{k}),  D_n(nu))=0$. Thus,
(3.12) still holds.
By Lemma 10, for any event ${\cal A}\in {\cal F}_{k}$,  $I(\Gamma_{k})I({\cal A})$ and
$T_\tau(S(\Gamma_{k}),  D_n(nu))$ are independent. Note
also that $I(\Gamma_k)$ is ${\cal F}_k$-measurable, so
\begin{eqnarray*}
&&{\bf E} (T_\tau(S(\Gamma_{k}),  D_n(nu))I(\Gamma_{k})\,\,|\,\, {\cal F}_{k})\\
&=&{\bf E} (T_\tau(S(\Gamma_{k}),  D_n(nu))({\bf E} I(\Gamma_{k})\,\,|\,\,{\cal F}_{k})=[{\bf E} T_\tau(S(\Gamma_{k}),  D_n(nu))]I(\Gamma_{k}). \hskip 1.2in (3.13)
\end{eqnarray*}
By (3.12) and (3.13),
\begin{eqnarray*}
&&\sum_{\Gamma_{k}} [{\bf E}T_\tau( S(\Gamma_{k}), D_n(nu))]I(\Gamma_{k} )\\
&\leq &\sum_{\Gamma_{k}}{\bf E} (T_\tau( \Gamma_{k},  D_n(nu))I(\Gamma_{k})\,\,| \,\,{\cal F}_{k})\\
&\leq & \sum_{\Gamma_{k}}[{\bf E}T_\tau( S(\Gamma_{k}), D_n(nu))]I(\Gamma_{k})+C\log^{2+2\delta}(n). \hskip 2.4in (3.14)
\end{eqnarray*}

Note that $B_\tau(k)\subset B_\tau(k+1)$, so by a subadditive property,
\begin{eqnarray*}
&& T_\tau(B_\tau(k+1),  D_n(nu))\leq T_\tau(B_\tau(k),  D_n(nu))\\
&\leq & T_\tau(B_\tau(k), B_\tau(k+1))+T_\tau(B_\tau(k+1),  D_n(nu)).\hskip 2.3in {(3.15)}
\end{eqnarray*}
By Lemma 5 and  Lemma 3,
$$T_\tau(B_\tau(k), B_\tau(k+1))\leq C\log^{2+2\delta}(n).\eqno{(3.16)}$$
Therefore, by (3.15)-(3.16) together with the same arguments of (3.11)-(3.14) for $B_\tau(k+1)$ and ${\cal F}_{k+1}$,
\begin{eqnarray*}
&&\sum_{\Gamma_{k+1}} [{\bf E}T_\tau( S(\Gamma_{k+1}), D_n(nu))]I(\Gamma_{k+1} )\\
&\leq &\sum_{\Gamma_{k+1}}{\bf E} (T_\tau( \Gamma_{k+1},  D_n(nu))I(\Gamma_{k+1})\,\,| \,\,{\cal F}_{k+1})\\
&\leq & \sum_{\Gamma_{k+1}}[{\bf E}T_\tau( S(\Gamma_{k+1}), D_n(nu))]I(\Gamma_{k+1})+C\log^{2+2\delta}(n). \hskip 2.1in (3.17)
\end{eqnarray*}
By Lemma 4, $\{B_\tau(k)=\Gamma_{k}\}$ are disjoint and $\cup_{\Gamma_k}\{B_\tau(k)=\Gamma_{k}\}=\Omega$. Thus,
by (3.17),
\begin{eqnarray*}
&&\sum_{\Gamma_{k+1}, \Gamma_{k}} [{\bf E}T_\tau( S(\Gamma_{k+1}), D_n(nu))]I(\Gamma_{k+1} )I(\Gamma_k)\\
&=&\sum_{\Gamma_{k+1}} [{\bf E}T_\tau( S(\Gamma_{k+1}), D_n(nu))]I(\Gamma_{k+1} )\\
&\leq &\sum_{\Gamma_{k+1}}{\bf E} (T_\tau( \Gamma_{k+1},  D_n(nu))I(\Gamma_{k+1})\,\,| \,\,{\cal F}_{k+1})\\
&\leq & \sum_{\Gamma_{k+1}}[{\bf E}T_\tau( S(\Gamma_{k+1}), D_n(nu))]I(\Gamma_{k+1})+\log^{2+2\delta}(n) \\
&=& \sum_{\Gamma_{k+1},\Gamma_k}[{\bf E}T_\tau( S(\Gamma_{k+1}), D_n(nu))]I(\Gamma_{k+1})I(\Gamma_k)+\log^{2+2\delta}(n), \hskip 1.7in (3.18)
\end{eqnarray*}
where the sums take over all possible $\Gamma_k$ and 
$\Gamma_{k+1}$ such that
$${\bf 0} \in \Gamma_k\subset \Gamma_{k+1} \mbox{ and }
d(\Gamma_k, \partial(\Gamma_{k+1})\leq \log^{1+\delta} (n).\eqno{(3.19)}$$
Similarly,
\begin{eqnarray*}
&&\sum_{\Gamma_{k+1}, \Gamma_{k}} [{\bf E}T_\tau( S(\Gamma_{k}), D_n(nu))]I(\Gamma_{k+1} )I(\Gamma_k)\\
&=&\sum_{\Gamma_{k}} [{\bf E}T_\tau( S(\Gamma_{k}), D_n(nu))]I(\Gamma_{k} )\\
&\leq &\sum_{\Gamma_{k}}{\bf E} (T_\tau(\Gamma_{k},  D_n(nu))I(\Gamma_{k})\,\,| \,\,{\cal F}_{k})\\
&\leq & \sum_{\Gamma_{k}}[{\bf E}T_\tau( S(\Gamma_{k}), D_n(nu))]I(\Gamma_{k})+\log^{2+2\delta}(n) \\
&=& \sum_{\Gamma_{k+1},\Gamma_k}[{\bf E}T_\tau( S(\Gamma_{k}), D_n(nu))]I(\Gamma_{k+1})I(\Gamma_k)+\log^{2+2\delta}(n).\hskip 1.8in (3.20)
\end{eqnarray*}
Note that
$$d(S(\Gamma_{k+1}), S(\Gamma_{k})) \leq C\log^{1+\delta} (n),$$
so by the same trick from (3.12), if $\Gamma_k$
and $\Gamma_{k+1}$ satisfy (3.19), we have
$$|{\bf E}T_\tau( S(\Gamma_{k+1}), D_n(nu))-{\bf E}T_\tau( S(\Gamma_{k}), D_n(nu))|\leq C\log^{2+2\delta}(n).\eqno{(3.21)}$$
If we substitute (3.10), (3.18), (3.20), and (3.21) into (3.9),  
we have 
\begin{eqnarray*}
&&|{\bf E}[T_\tau(D_n({\bf 0}), D_n(nu))I_2(k)\,\, |\,\, {\cal F}_{k+1}]
 -{\bf E}[T_\tau(D_n({\bf 0}), D_n(nu))I_2(k)\,\,|\,\,{\cal F}_k]|\\
&= & |{\bf E} [T_\tau( B_\tau(k), D_n(nu) )\,\,|\,\,{\cal F}_{k+1}]-
{\bf E} [T_\tau(B_\tau(k), D_n(nu))\,\,| \,\,{\cal F}_k]|\\
&\leq & \sum_{\Gamma_{k+1}, \Gamma_k}\left(|{\bf E}T_\tau( S(\Gamma_{k+1}), D_n(nu))-{\bf E}T_\tau( S(\Gamma_{k}), D_n(nu))| I(\Gamma_{k+1})I(\Gamma_k)\right)+C_1\log^{2+2\delta}(n)\\
&\leq & C_2 \log^{2+2\delta}(n).\hskip 4.8in (3.22)
\end{eqnarray*}

By (3.2), (3.3), (3.6), and (3.22), together with Azuma's lemma, for all $x >0$,
\begin{eqnarray*}
&&{\bf P} \left[  |T_\tau(D_n({\bf 0}), D_n(nu)) -{\bf E}T_\tau(D_n({\bf 0}), D_n(nu))  |\geq x\sqrt{n}\log^{3+3\delta} n\right] \\
&\leq & C_1 \exp \left( -2 x^2{n \log^{6+6\delta}n \over (C\log^{2+2\delta} n )^2 n \log^{1+\delta} n)}\right)\\
&\leq &C_1 \exp\left(-C_2x^2\log^{1+\delta}n\right).\hskip 4in (3.23)
\end{eqnarray*}
We use (3.23) to show Theorems 1 and 2.\\

\section{ Proofs of Theorems 1 and 2.}
By Lemma 8, 
$${\bf P}[ {T}(D_n({\bf 0}), D_n(nu)) \neq T_\tau(D_n({\bf 0}), D_n(nu))]\leq C_1 \exp\left(-C_2\log^{1+\delta} n\right).\eqno{(4.1)}$$
Let ${\cal D}_n$ be the event that $T(D_n({\bf 0}), D_n(nu)) \neq T_\tau(D_n({\bf 0}), D_n(nu))$.
\begin{eqnarray*}
{\bf  E}(T_\tau(D_n({\bf 0}), D_n(nu)))&=&{\bf E}(T_\tau(D_n({\bf 0}), D_n(nu)))I({\cal D}_n)+{\bf E}(T_\tau(D_n({\bf 0}), D_n(nu)))[1-I({\cal D}_n)]\\
&\leq& {\bf E}(T(D_n({\bf 0}), D_n(nu)))+{\bf E}(T_\tau({\bf 0}, nu))[1-I({\cal D}_n)].
\end{eqnarray*}
Now we estimate ${\bf E}(T_\tau(D_n({\bf 0}), D_n(nu)))[1-I({\cal D}_n)]$. By Lemma 6, we have
$$T_\tau(D_n({\bf 0}), D_n(nu))\leq C n \log^{1+\delta}(n).\eqno{(4.2)}$$
By (4.1) and (4.2),
\begin{eqnarray*}
&&{\bf E}(T_\tau(D_n({\bf 0}), D_n(nu)))[1-I({\cal D}_n)]\\
&\leq & \left[Cn \log^{1+\delta}(n)\right]\left[C_1 \exp(-C_2\log^{1+\delta} n)\right]\\
&\leq & C_3 \exp\left(C_4 \log^{1+\delta} (n)\right).\hskip 4in {(4.3)}
\end{eqnarray*}
Therefore, we have
$${\bf  E}T_\tau(D_n({\bf 0}), D_n(nu))
\leq {\bf E}T(D_n({\bf 0}), D_n(nu))+C_1\exp\left(-C_2 \log^{1+\delta} (n)\right).\eqno{(4.4)}
$$
Similarly,
$${\bf  E}T(D_n({\bf 0}), D_n(nu))\leq  {\bf E}T_\tau(D_n({\bf 0}), D_n(nu))+{\bf E}T(D_n({\bf 0}), D_n(nu))
[1-I({\cal D}_n)].\eqno{(4.5)}$$
It is more difficult to estimate ${\bf E}T(D_n({\bf 0}), D_n(nu))[1-I({\cal D}_n)]$ since we do not have a bound in (4.2)
for $T(D_n({\bf 0}), D_n(nu)).$
Note that  $nu\in [-n,n]^d$, so
$${\bf E}T(D_n({\bf 0}), D_n(nu))[1-I({\cal D}_n)] \leq {\bf E} \left[ \sum_{e\in [-n,n]^d} t(e) \right][1-I({\cal D}_n)].\eqno{(4.6)}$$
Let ${\cal G}(e)$ be the event that
$$ t(e) \geq C_1 \exp\left((C_2 \log^{1+\delta} n)/2\right),$$
where $C_1$ and $C_2$ are in (4.1).
With this definition, by (4.1),
\begin{eqnarray*}
&&{\bf E} \left[ \sum_{e\in [-n,n]^d} t(e) \right][1-I({\cal D}_n)]\\
&\leq & n^{2d} C_1 \exp\left((C_2 \log^{1+\delta} n)/2\right)C_1 \exp(-C_2 \log^{1+\delta} n)+ {\bf E} \left[ \sum_{e\in [-n,n]^d} t(e) I({\cal G}(e))\right].\hskip 0.5in {(4.7)}
\end{eqnarray*}
Let us estimate ${\bf E} \left[ \sum_{e\in [-n,n]^d} t(e) I({\cal G}_n)\right]$:
$${\bf E} \left[ \sum_{e\in [-n,n]^d} t(e) I({\cal G}(e))\right]\leq n^{2d} {\bf E} [ t(e) ]I\left(t(e)  \geq C_1 \exp\left({C_2 \over 2}\log^{1+\delta} n\right)\right).\eqno{(4.8)}$$
Note that ${\bf E}(t(e))^{1+\eta} <\infty$, so by Markov's inequality,
\begin{eqnarray*}
&&{\bf E} [ t(e) ]I\left(t(e)  \geq C_1 \exp\left({C_2 \over 2}\log^{1+\delta} n\right)\right)\\
&\leq &\int_{t\geq C_1 \exp\left({C_2 \over 2} \log^{1+\delta} n\right)}
{\bf P} [ t(e) \geq t] dt \\
&\leq &C_3\int_{t\geq C_1 \exp\left({C_2 \over 2} \log^{1+\delta} n\right)}t^{-1-\eta} dt \\
&\leq & C_4 \exp\left(-\eta (C_2 \log^{1+\delta} n)/2\right).\hskip 9cm (4.9)
\end{eqnarray*}
If we substitute (4.9) into (4.8), we have
$${\bf E} \left[ \sum_{e\in [-n,n]^d} t(e) I({\cal G}(e))\right]\leq C_1 \exp\left(-\eta (C_2 \log^{1+\delta} n)/2\right).\eqno{(4.10)}$$
By (4.6) and (4.10), we have
$${\bf  E}(T(D_n({\bf 0}), D_n(nu)))\leq  {\bf E}(T_\tau(D_n({\bf 0}), D_n(nu)))+C_1 \exp\left(-\eta C_2 \log^{1+\delta} (n)\right).\eqno{(4.11)}$$
Together with (4.2),  we have
$$|{\bf E}T(D_n({\bf 0}), D_n(nu)) - {\bf E}T_\tau(D_n({\bf 0}), D_n(nu))|\leq C_1 \exp\left(-C_2 \log^{1+\delta} n\right).\eqno{(4.12)}$$

By (4.12) and (3.23), again,  for all large $n$,
\begin{eqnarray*}
&&{\bf P}\left [ \left |T(D_n({\bf 0}), D_n(nu)) -{\bf E}T(D_n({\bf 0}), D_n(nu)) \right|\geq  \sqrt{n} \log^{3+3\delta} n\right] \\
&\leq &{\bf P}\left [ |T(D_n({\bf 0}), D_n(nu)) -{\bf E}T_\tau(D_n({\bf 0}), D_n(nu))|\geq  {\sqrt{n}\over 2}\log^{3+\delta}n\right] \\
&\leq &{\bf P}\left [ |T(D_n({\bf 0}), D_n(nu)) -{\bf E}T_\tau(D_n({\bf 0}), D_n(nu))|\geq  {\sqrt{n}\over 2}\log^{3+\delta}n, {\cal D}_n\right] 
+C_1\exp\left(-C_2  \log^{1+\delta} n\right)\\
&= &{\bf P}\left [ |T_\tau(D_n({\bf 0}), D_n(nu)) -{\bf E}T_\tau(D_n({\bf 0}), D_n(nu))|\geq  {\sqrt{n}\over 2}\log^{3+\delta}n, {\cal D}_n\right] 
+C_1\exp\left(-C_2 \log^{1+\delta} (n)\right)\\
&\leq & C_3 \exp\left(-C_4 \log^{1+\delta}(n)\right).\hskip 4in (4.13)
\end{eqnarray*}

Note that if 
$$|T(D_n({\bf 0}), D_n(nu)) - {\bf E}T(D_n({\bf 0}), D_n(nu))|^{2m}\geq k\eqno{(4.14)}$$
for $m\geq 1$ and
$$k \geq 2 n^{4dm} \exp\left(m \log^{1+\delta} n\right),\eqno{(4.15)}$$
then  by (4.2) and (4.12),
$$T(D_n({\bf 0}), D_n(nu)) \geq k^{1/2m}/2.$$
Note also that
$$k^{1/2m}/2\leq 
T(D_n({\bf 0}), D_n(nu)) \leq \sum_{e\in [-n,n]^d} t(e),\eqno{(4.16)}$$
so by (4.16), if (4.15) holds, then there exists an edge $e\in
[-n,n]^d$ such that
$$t(e) \geq {k^{1/2m}\over 2n^{2d}}.\eqno{(4.17)}$$
Thus,
\begin{eqnarray*}
&&\sum_{k \geq  (2n)^{4dm} \exp\left(m \log^{1+\delta} n\right)} {\bf P} [|T(D_n({\bf 0}), D_n(nu)) - 
{\bf E}T(D_n({\bf 0}), D_n(nu))|^{2m}\geq k]\\
&\leq & n^{2d}\sum_{k^{1/2m} \geq  2n^{2d} \exp\left(\log^{1+\delta} n\right)} {\bf P} [t(e)\geq {k^{1/2m}\over 2n^{2d}}]\\
&\leq &  n^{2d}\sum_{k^{1/2m} /2n^{2d}\geq  \exp\left(\log^{1+\delta} n\right)} {\bf P} [t(e)\geq {k\over 2n^{2d}}]\\
&\leq &C n^{2d}{\bf E} t(e) I\left[ t(e)\geq  \exp\left( \log^{1+\delta} n\right)\right].
\end{eqnarray*}
By using  (4.9) above, we have
\begin{eqnarray*}
&&\sum_{k \geq (2 n)^{2dm} \exp\left(m \log^{1+\delta} n\right)} {\bf P} [|T(D_n({\bf 0}), D_n(nu)) - 
{\bf E}T(D_n({\bf 0}), D_n(nu))|^{2m}\geq k]\\
&\leq &
C_1\exp\left(-{C_2} \log^{1+\delta}(n)\right)\hskip 4in {(4.18)}
\end{eqnarray*}

{\bf Proof of Theorem 2.}
The probability estimate in Theorem 2 follows from (4.13) directly. Now we show the mean estimate in Theorem 2.
Note that
\begin{eqnarray*}
&&{\bf E} \left[T(D_n({\bf 0}), D_n(nu)) -{\bf E}T(D_n({\bf 0}), D_n(nu))\right]^{2m}\\
&\leq & C\sum_{k=1}^\infty {\bf P} \left[ (T(D_n({\bf 0}), D_n(nu)) -{\bf E}T(D_n({\bf 0}), D_n(nu)))^{2m}\geq k\right]\\
&=& C\sum_{k=1}^{n^m \log^{7m} n}  {\bf P} \left[(T(D_n({\bf 0}), D_n(nu)) -{\bf E}T(D_n({\bf 0}), D_n(nu)))^{2m}\geq k\right]\\
& +&C\sum_{n^m \log^{7m} n}^{(2 n)^{2dm} \exp\left(m \log^{1+\delta} n\right)} {\bf P} \left[(T(D_n({\bf 0}), D_n(nu)) -{\bf E}T(D_n({\bf 0}), D_n(nu)))^{2m}\geq k\right]\\
& +&C\sum_{(2 n)^{2dm} \exp\left(m\log^{1+\delta} n\right)}^\infty {\bf P} \left[(T(D_n({\bf 0}), D_n(nu)) -{\bf E}T(D_n({\bf 0}), D_n(nu)))^{2m}\geq k\right]\\
&=& I+II+III,
\end{eqnarray*}
where $I$, $II$ and $III$ are the first,  second and  third
sums, respectively. Clearly,
$$I \leq  C n^m \log^{7m} (n).$$
By (4.13), if we  take $6\delta \leq 0.5$, we have
$$II \leq (2 n)^{2dm} C_1\exp\left(m\log^{1+\delta} n\right) C_3\exp\left(-{C_4} \log^{1+2\delta} n\right)\leq
C_5 \exp(-C_6  \log^{1+2\delta} n).$$
By (4.18),
$$III \leq C_1\exp\left(-{C_2} \log^{1+\delta}(n)\right).$$
Thus, there exists $C=C(F, d, m, M, \eta)$ such that
$${\bf E} \left[T(D_n({\bf 0}), D_n(nu)) -{\bf E}T(D_n({\bf 0}), D_n(nu))\right]^{2m}\leq C n^m \log^{7m} (n).\eqno{(4.19)}$$
Therefore, the mean estimate in Theorem 2 follows
from (4.19). $\Box$\\

Now we show Theorem 1.  By the definition,
$$  T(D_n({\bf 0}), D_n(nu))\leq T({\bf 0}, nu).\eqno{(4.20)}$$
Therefore,
$${\bf E} T(D_n({\bf 0}), D_n(nu)) \leq {\bf E} T({\bf 0}, nu).\eqno{(4.21)}$$
On the other hand, suppose that there exists an optimal path of $ T(D_n({\bf 0}), D_n(nu)) $ from $v$ at $\partial D_n({\bf 0})$ to 
$v'$ at $\partial D_n(nu)$. We find  deterministic  paths $\gamma$ and $\gamma'$ from the origin to $v$ and from $nu$ to $v'$ with less than
$2 \log^{1+\delta} n$ edges.
Thus,
$$T({\bf 0}, nu)\leq T(D_n({\bf 0}), D_n(nu))+ T(\gamma)+T(\gamma').\eqno{(4.22)}$$
As we defined,
$${\bf E} T(\gamma')+T(\gamma')\leq  2(\log^{1+\delta} n ){\bf E} t(e) .\eqno{(4.22)}$$

{\bf Proof of Theorem 1.}  By (4.20) and (4.22),
\begin{eqnarray*}
&&{\bf E} \left[T({\bf 0}, nu) -{\bf E}{T}({\bf 0}, nu)\right]^{2m}\\
&\leq & {\bf E} [T({\bf 0}, nu) -T(D_n({\bf 0}), D_n(nu))+ T(D_n({\bf 0}), D_n(nu)) -{\bf E} T(D_n({\bf 0}), D_n(nu))\\
&&+{\bf E} T(D_n({\bf 0}), D_n(nu))  -{\bf E}{T}({\bf 0}, nu)]^{2m}\\
&\leq & C\log^{2m(1+\delta)}(n) {\bf E}(T(D_n({\bf 0}), D_n(nu)) -{\bf E} T(D_n({\bf 0}), D_n(nu)))^{2m}+ C\log^{1+\delta} (n).
\end{eqnarray*}
Therefore, Theorem 1 follows from (4.19) by selecting small $\delta$. $\Box$

\section {Proof of Theorem 3.}
We follow the idea in Zhang (2005) to show Theorem 3.
For each $n >0$, we define face-face first passage time  as
$$\Phi_{k,m}(n)=\inf\{T(\gamma): \gamma \mbox{ is a path from $\{k\}\times [-n^2,n^2]^{d-1}$ to $\{m\} \times [-n^2,n^2]^{d-1}$ }\},$$
where we require that all paths stay inside $(k,m)\times [-n^2, n^2]^{d-1}$ except their ending points.
Here without loss of generality, we assume that  $\log^{1+\delta} n$ and $n/ \log^{1+\delta} n$ are integers.
Now we divide (see Fig. 3) $ [-n,n]^{d-1}$ into $(d-1)$-dimensional cubes with 
a side length of $\log^{1+\delta} n$. We denote them by
$S_1,\cdots, S_l$
for  
$$l =O( (n^2/\log^{1+\delta} n)^{d-1}).\eqno{(5.1)}$$
Let
$$S_i(k)=\{k\}\times S_i.$$
In particular, let $S_0(0)$ (see Fig. 3) be the cube that
contains the origin.
Note that there always exists a path such that
$$T(\gamma)=\Phi_{0,n}(n)\mbox{ for each configuration}.\eqno{(5.2)}$$
For $v_n\in \{0\}\times [-n^2, n^2]$
 and $v_n'\in \{n\} \times [-n^2,n^2]^{d-1}$,
let us set the following event:
$${\cal A}(v_n,v'_n)=\{\, \exists\,\, \mbox{ a path } \gamma\subset (0, n)\times [-n^2, n^2]
\mbox{ from $v_n$ to $v_n'$ such that }T(\gamma)= \Phi_{0,n}(n)\}.$$
It follows from (5.2) that
$${\bf P}\left[\bigcup_{u_n,u_n'} {\cal A} (u_n, u_n')\right]=1.$$
Note that there are at most $Cn^{2(d-1)}$ choices for $u_n$ and $C n^{2(d-1)}$ choices for $u'_n$,
 so there exist fixed $v_n$ and $v_n'$ such that
$${\bf P}\left[{\cal A}(v_n,v_n')\right]\geq {C\over n^{4d}}.\eqno{(5.3)}$$
We assume that $v_n\in S_j(0)$ for a fixed $j$ and
$${\cal A}(S_j(0),v_n')={\cal A}(v_n,v_n').$$
We also set
\begin{eqnarray*}
&&{\cal A}'(v_n', S_j(2n))\\
&=&\mbox{\{the reflection of  all configurations in the event ${\cal A}(S_j(0), v_n')$ about the hyperplane $x_1=n$\}}.
\end{eqnarray*}

Note that  two events ${\cal A}(S_j(0),v'_n)$ and $ {\cal A}'(v_n', S_j(2n))$ are independent since they use the edges 
in different vertex sets.
Note also that on $ {\cal A}'(v_n', S_j(2n))$, there exists a path $\gamma'$ from $v_n'$
to $S_j(2n)$ inside $(n,2n)\times [-n^2,n^2]^{d-1}$, except at its two starting points
with
$$T(\gamma')= \Phi_{n,2n}(n).$$

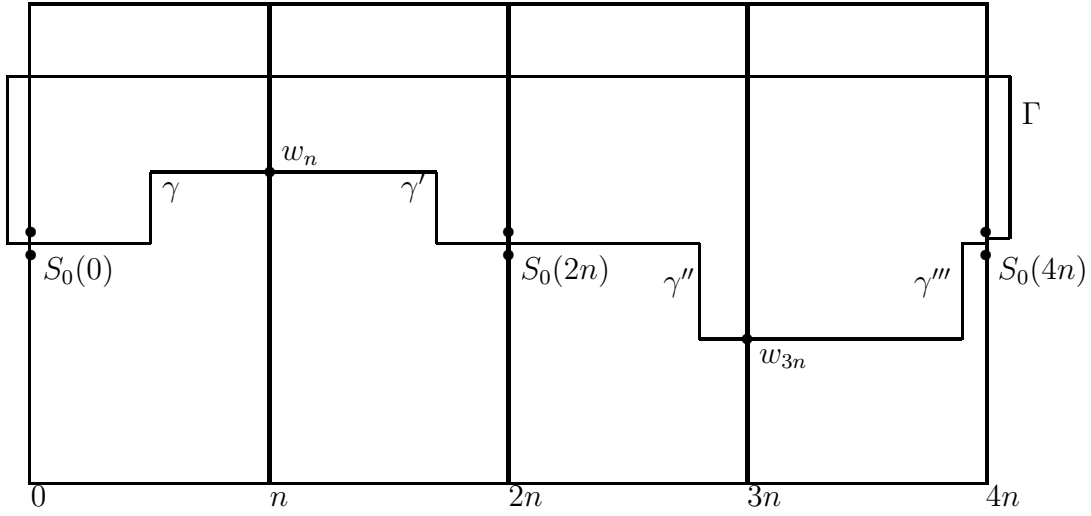
\begin{figure}\label{F:graphG}
\begin{center}
\setlength{\unitlength}{0.0125in}%
\begin{picture}(200,180)(67,910)
\thicklines
\put(150, 900){\framebox(100,200)[br]{\mbox{$$}}}
\put(50, 900){\framebox(100,200)[br]{\mbox{$$}}}
\put(-50, 900){\framebox(100,200)[br]{\mbox{$$}}}
\put(250, 900){\framebox(100,200)[br]{\mbox{$$}}}

\put(150,995){\circle*{4}}
\put(150,1005){\circle*{4}}
\put(155,985){$S_0(2n)$}

\put(-50,995){\circle*{4}}
\put(-50,1005){\circle*{4}}
\put(-45,985){$S_0(0)$}

\put(350,995){\circle*{4}}
\put(350,1005){\circle*{4}}
\put(355,985){$S_0(4n)$}

\put(-50,1000){\line(1,0){50}}
\put(0,1000){\line(0,1){30}}
\put(0,1030){\line(1,0){50}}
\put(50,1030){\circle*{4}}
\put(55,1035){$w_n$}
\put(5,1020){$\gamma$}

\put(50,1030){\line(1,0){70}}
\put(120,1030){\line(0,-1){30}}
\put(120,1000){\line(1,0){30}}
\put(105,1020){$\gamma'$}

\put(150,1000){\line(1,0){80}}
\put(230,1000){\line(0,-1){40}}
\put(230,960){\line(1,0){20}}
\put(255,950){$w_{3n}$}
\put(250, 960){\circle*{4}}
\put(215,980){$\gamma''$}

\put(250,960){\line(1,0){90}}
\put(340,960){\line(0,1){40}}
\put(340,1000){\line(1,0){10}}
\put(320,980){$\gamma'''$}

\put(-50,1000){\line(-1,0){10}}
\put(-60,1000){\line(0,1){70}}
\put(-60,1070){\line(1,0){420}}
\put(360,1070){\line(0,-1){68}}
\put(360,1002){\line(-1,0){10}}
\put(365,1050){$\Gamma$}

\put(-50,890){$0$}
\put(50,890){$n$}
\put(150,890){$2n$}
\put(250,890){$3n$}
\put(350,890){$4n$}

\end{picture}
\end{center}
\caption{\em $\gamma$, $\gamma'$, $\gamma''$, and $\gamma'''$ are optimal paths of $\Phi_{0,n}$, $\Phi_{n, 2n}$, 
$\Phi_{2n, 3n}$, and $\Phi_{3n, 4n}$, respectively. On ${\cal A}_1\cap {\cal A}_2\cap {\cal A}_3\cap {\cal A}_4$,
$\gamma\cup \gamma'$ consists of a path
 $s_{0, 2n}$, and $\gamma''\cup \gamma'''$ consists of a path  $s_{2n, 4n}$. On ${\cal B}_n$,
$\Gamma$ is an optimal path
 $s_{0, 4n}$ that crosses the four boxes in the first coordinate direction. }
\end{figure} 

By translation invariance, 
with a probability larger than ${C\over n^{8d}}$, there exist $\gamma$ from $S_0(0)$ to $w_n$ and
$\gamma'$ from $w_n$ to $S_0(n)$ inside $(0,n)\times [-n^2, n^2]^{d-1}$ and $(n,2n)\times [-n^2, n^2]^{d-1}$
(see Fig. 3),
respectively, such that
$$T(\gamma)=\Phi_{0,n}(n) \mbox{ and } T(\gamma')=\Phi_{n,2n}(n).\eqno{(5.4)}$$
We denote the events that there exist 
$\gamma$ and $\gamma'$ satisfying (5.4) by ${\cal A}(S_0(0),w_n)$ and $ {\cal A}'(w_n,S_0(2n))$, respectively.
Let 
$$s_{0,n}=\inf_{\gamma} \{T(\gamma): \gamma \mbox{ is a path from $S_0(0)$ to $S_0(n)$}\}.$$
With this definition, using the exact same proof of Theorem 2
and (4.18), we have the following lemma:\\

{\bf Lemma 11.} {\em If $F(0) < p_c$ and ${\bf E}(t(e))^{1+\eta} < \infty$ for  $\eta >0$, then 
there exists $C_i=C_i(d, \eta, \delta, F)$ for $i=1,2$ such that
$${\bf P}\left [ \left |{s}_{0,n} -{\bf E}s_{0,n} 
\right|\geq \sqrt{n} \log^{4} n \right] \leq C_1\exp\left(-C_2 \log^{1+\delta} (n)\right).$$
Furthermore,}
\begin{eqnarray*}
&&\sum_{k \geq (2 n)^{2d} \exp\left( \log^{1+\delta} n\right)} {\bf P} [|s_{0,n} - {\bf E}s_{0, n}|\geq k]\\
&\leq &
C_1\exp\left(-{C_2} \log^{1+\delta}(n)\right)
\end{eqnarray*}

Note that $\gamma\,\,\cup \,\,\gamma'$ is a path from $S_0(0)$ to $S_0(2n)$. Therefore, on 
${\cal A}(S_0(0),w_n)\,\,\cap\,\, {\cal A}'(w_n,S_0(2n))$,
$$s_{0,2n}\leq T(\gamma\cup \gamma')=\Phi_{0, n}(n)+ \Phi_{n,2n}(n).\eqno{(5.5)}$$
We may continue to work with this method on the cylinders $(2n, 3n)\times [-n^2, n^2]^{d-1}$ and $(3n, 4n)\times [-n^2, n^2]^{d-1}$. 
With the same argument, we can show that (see Fig. 3)
on ${\cal A}(S_0(0), w_n)\cap {\cal A}'( w_n,S_0(2n))\cap  {\cal A}(S_{0}(2n), w_{3n})\cap {\cal A}'( w_{3n},S_{0}(4n))$,
there exist $\gamma$ from $S_0(0)$ to $w_n$, and
$\gamma'$ from $w_n$ to $S_0(n)$
inside $(0,n)\times [-n^2, n^2]^{d-1}$, and $(n,2n)\times [-n^2, n^2]^{d-1}$, respectively. 
There also exist $\gamma''$ from $S_0(2n)$ to $w_{3n}$, and $\gamma'''$ from
$w_{3n}$ to $S_{0}(4n)$ inside $(2n ,3n)\times [-n^2, n^2]^{d-1}$, and $(3n, 4n)\times [-n^2, n^2]^{d-1}$,
respectively, such that
$$T(\gamma')=\Phi_{0,n}(n) \mbox{ and } T(\gamma')=\Phi_{n,2n}(n) \mbox{ and } T(\gamma'')=\Phi_{2n,3n}(n)\mbox{ and } T(\gamma''')=\Phi_{3n,4n}.$$
By this definition, on ${\cal A}(S_0(0), w_n)\cap {\cal A}'( w_n,S_0(2n))\cap  {\cal A}(S_{0}(2n), w_{3n})\cap {\cal A}'( w_{3n},S_{0}(4n))$,
$$s_{0,2n}+s_{2n, 4n}\leq T(\gamma\cup \gamma'\cup\gamma''\cup \gamma''')=\Phi_{0,n}(n)+\Phi_{n,2n}(n)+\Phi_{2n,3n}(n)+\Phi_{3n,4n}(n)
.\eqno{(5.6)}$$
By independence, 
$${\bf P} \left[{\cal A}(S_0(0), w_n)\cap {\cal A}'( w_n,S_0(2n))\cap  {\cal A}(S_{0}(2n), w_{3n})\cap {\cal A}'( w_{3n},S_{0}(4n))\right] \geq {C\over n^{16d}}.\eqno{(5.7)}$$
We now consider an optimal path $\gamma$ for $s_{0, 4n}$. By Proposition 5.8 of Kesten (1986), 
$${\bf P}\left[\exists\,\,\,\gamma \in (-\infty, \infty) \times [-n^2, n^2]^{d-1}\mbox{ with }
T(\gamma)= s_{0, 4n}\right ]\geq 1-C_1\exp(-C_2 n).\eqno{(5.8)}$$
We denote by ${\cal B}_n$ the event in (5.8).
On ${\cal B}_n$, note that any optimal path of $s_{0, 4n}$ must cross out, in the first coordinate direction,
 the box $[0, 4n]\times [-n^2, n^2]^d$, so
$$\Phi_{0,n}+\Phi_{n,2n}+\Phi_{2n,3n}+\Phi_{3n,4n}\leq s_{0, 4n}.\eqno{(5.9)}$$
For simplicity, we write,
$${\cal A}(S_0(0), w_n)\cap {\cal A}'( w_n,S_0(2n))\cap  {\cal A}(S_{0}(2n), w_{3n})\cap {\cal A}'( w_{3n},S_{0}(4n))={\cal A}_1\cap {\cal A}_2
\cap {\cal A}_3\cap {\cal A}_4.$$
By (5.6) and (5.9), 
\begin{eqnarray*}
&&{\bf E} \left(s_{0,2n}\,\,|\,\,{\cal A}_1\cap {\cal A}_2\cap {\cal A}_3\cap {\cal A}_4 \cap {\cal B}_n\right)
+{\bf E} \left(s_{2n, 4n}\,\,|\,\,{\cal A}_1\cap {\cal A}_2\cap {\cal A}_3\cap {\cal A}_4 \cap {\cal B}_n\right)\\
&\leq &{\bf E} \left(s_{0,4n}\,\,|\,\,{\cal A}_1\cap {\cal A}_2\cap {\cal A}_3\cap {\cal A}_4 \cap {\cal B}_n\right).\hskip 2.9in (5.10)
\end{eqnarray*}

Now we estimate 
\begin{eqnarray*}
&&{\bf E}\left (|{\bf E}s_{0, 2n}-s_{0, 2n}|\,\,\, |\,\,\,{\cal A}_1\cap {\cal A}_2\cap {\cal A}_3\cap {\cal A}_4 \cap {\cal B}_n\right)\\
&\leq & 
C\sum_{m\geq \sqrt{n} (\log^{1+\delta}  n)}
{{\bf P} \left[|{\bf E}s_{0, 2n}-s_{0, 2n}|\,\,\geq m\right]\over {\bf P}\left[{\cal A}_1\cap {\cal A}_2\cap {\cal A}_3\cap {\cal A}_4 \cap {\cal B}_n\right]}\\
&+& C\sum_{m<  \sqrt{n} (\log^{4}  n)}
{{\bf P} \left[|{\bf E}s_{0, 2n}-s_{0, 2n}|\,\,\geq m\,\,|\,\,{\cal A}_1\cap {\cal A}_2\cap {\cal A}_3\cap {\cal A}_4 \cap {\cal B}_n\right]}\\
&\leq & 
C\sum_{m\geq \sqrt{n} (\log^{4}  n)}
{{\bf P} \left[|{\bf E}s_{0, 2n}-s_{0, 2n}|\,\,\geq m\right]\over {\bf P}\left[{\cal A}_1\cap {\cal A}_2\cap {\cal A}_3\cap {\cal A}_4 \cap {\cal B}_n\right]}+ C\sqrt{n} \log^{4}(n). \hskip 4cm (5.11)
\end{eqnarray*}
By  (5.7) and (5.8), 
$${\bf P}^{-1}\left[{\cal A}_1\cap {\cal A}_2\cap {\cal A}_3\cap {\cal A}_4 \cap {\cal B}_n\right] \leq Cn^{16d}.\eqno{(5.12)}$$
By Lemma 11, 
$$\sum_{m \geq \sqrt{n} (\log^{1+\delta}n)}{\bf P} \left[|{\bf E}s_{0, 2n}-s_{0, 2n}|\geq m\right ]\leq C_1\exp(-C_2  \log^{1+\delta}(m)).\eqno{(5.13)}$$
We substitute (5.12) and (5.13) into (5.11) to show
$${\bf E} \left(|{\bf E}s_{0,2n}- s_{0,2n}|\,\,|\,\,{\cal A}_1\cap {\cal A}_2\cap {\cal A}_3\cap {\cal A}_4 \cap {\cal B}_n\right)\leq 
C\sqrt{n} \log^{4}(n).\eqno{(5.14)}$$

By (5.14), 
\begin{eqnarray*}
&&|{\bf E}s_{0, 2n}-{\bf E}\left[s_{0, 2n}\,\,|\,\, {\cal A}_1\cap {\cal A}_2\cap {\cal A}_3\cap {\cal A}_4 \cap {\cal B}_n\right]|\\
& =& |{\bf E}[({\bf E} s_{0, 2n}-s_{0, 2n})\,\,|\,\, 
{\cal A}_1\cap {\cal A}_2\cap {\cal A}_3\cap {\cal A}_4 \cap {\cal B}_n]|\\
&\leq & {\bf E}\left( |{\bf E}s_{0,2n}- s_{0,2n}|\,\,|\,\,{\cal A}_1\cap {\cal A}_2\cap {\cal A}_3\cap {\cal A}_4 \cap {\cal B}_n\right)\\
&\leq &
C\sqrt{n} \log^{4}(n). \hskip 4.5in (5.15)
\end{eqnarray*}
We use the same argument of (5.15) to show that
\begin{eqnarray*}
&&|{\bf E}s_{0, 4n}-{\bf E}\left[s_{0, 4n}\,\,|\,\, {\cal A}_1\cap {\cal A}_2\cap {\cal A}_3\cap {\cal A}_4 \cap {\cal B}_n\right]|\\
&\leq & {\bf E} (|{\bf E}s_{0,4n}- s_{0,4n}|\,\,|\,\,{\cal A}_1\cap {\cal A}_2\cap {\cal A}_3\cap {\cal A}_4 \cap {\cal B}_n)\\
&\leq &
C\sqrt{2n} \log^{4}(n). \hskip 4.5in (5.16)
\end{eqnarray*}

By (5.15), (5.16), and (5.10), we have
$${{\bf E }s_{0, 2n}\over 2n}\leq {{\bf E}s_{0, 4n}\over 4n} + {C\log^{4}(n)\over \sqrt{2n}}.\eqno{(5.17)}$$
If we iterate (5.17),  note that
$$\log 2^j+ \log^4 n \leq (\log 2^j )( \log^4 n) \mbox{ for  large $n$},$$
so
$${\bf E }{s_{0, n}\over n}\leq {{\bf E}s_{0, 2^i n}\over 2^i n}  +{C\log^{4} n\over \sqrt{n} }\sum_{j=0}^i {\log 2^j\over \sqrt{2^j}}.\eqno{(5.18)}$$

{\bf Proof of Theorem 3.}
By the same argument of (1.11),
$${\bf E} s_{0,n}\leq {\bf E} a_{0,n}\leq {\bf E} s_{0,n}+C\log^{4} (n).\eqno{(5.19)}$$
By (5.19),
$$
\lim_{n\rightarrow \infty} {\bf E} {s_{0,n}\over n}
=\mu_F \mbox{ a.s. and in } L_1.\eqno{(5.20)}
$$
If we let $i\rightarrow \infty$ in (5.18), by (5.19)
and  (5.20), we have
$${\bf E }{a_{0, n}\over n}\leq {\bf E }{s_{0, n}\over n}\leq \mu_F  +{C\log^{4} n\over \sqrt{n} }.\eqno{(5.21)}$$
On the other hand, by (1.2),
$$\mu_F \leq {\bf E }{a_{0, n}\over n}. \eqno{(5.22)}$$
Therefore, Theorem 3 follows from (5.21) and (5.22). $\Box$

\begin{center}
{\bf \large References}
\end{center}
Alexander, K. (1993). A  note on some rates of convergence in first-passage percolation. {\em Ann. Apple. Pro Bab.} {\bf 3}
81--91.\\
Benjamini, I., Kalai, G. and  Schramm, O. (2003). First passage percolation has sub linear distance variance. 
{\em Ann. Pro Bab.} {\bf 31}  1970--1978.\\
Cox, J. T. and Durrett, R. (1981). Some limit theorems for percolation processes with necessary and sufficient conditions. {\em Ann. Pro Bab.} {\bf 9} 583--603. \\
Grimmett, G. (1999). {\em Percolation.} Springer, Berlin.\\
 Hammersley, J. M. and  Welsh, D. J. A. (1965).
First-passage percolation, subadditive processes,
stochastic networks and generalized renewal theory.
In {\em Bernoulli, Bayes, Laplace Anniversary Volume} 
(J. Neyman   and L. LeCam  eds.) 61--110 Springer, Berlin.\\
Higuchi, Y. and Zhang, Y. (2000). On the speed of convergence for two dimensional first-passage Ising percolation.
{\em Ann. of Probab.} {\bf 28} 353--376.\\
Kesten, H. (1986). Aspects of first-passage percolation. {\em Lecture Notes in
Math.} {\bf 1180} 125--264. Springer, Berlin.\\
Kesten, H. (1993). On the speed of convergence in first passage percolation. {\em Ann Appl. Probab.} {\bf 3} 296--338.\\
Kesten, H. and Zhang, Y. (1997). A  central limit theorem  for critical first passage percolation in two dimensions. {\em PTRF} {\bf 107} 137--160.\\
Smythe, R. T. and Wierman, J. C. (1978).
First passage percolation on the square lattice.
{\em Lecture Notes in Math.} {\bf 671} Springer, Berlin.\\
Talagrand, M. (1995). Concentration of measure and isoperimetric inequalities in product spaces. {\em Inst. Hautes Publ. Math. Etudes Sci.}  {\bf 81} 73--205.\\
Zhang, Y. (1995). Supercritical behaviors in first-passage percolation. 
{\em Stoch. Proc. Appl.} {\bf 59} 251--266.\\
Zhang, Y. (2005). On the speeds of convergence and concentration of a subadditive
ergodic process. Preprint.\\
Zhang, Y. (2006). The divergence of fluctuations for  shape in first passage percolation. {\em Probab. Theory and Relat. Fields}  {\bf 136} 298--320.\\

\noindent
Yu Zhang\\
Department of Mathematics\\
University of Colorado\\
Colorado Springs, CO 80933\\
email: yzhang3@uccs.edu\\
\end{document}